\documentclass{article}

\usepackage{amsmath,amssymb}

\newtheorem{definition}{Definition}

\newtheorem{lemma}{Lemma}
\newtheorem{proposition}{Proposition}
\newtheorem{theorem}{Theorem}

\begin{document}

\title{Vanishing of the Logarithmic Trace \\of  Generalized Szeg\"o Projectors}

\author{L. Boutet de Monvel}
\date{}

\maketitle

\small{\noindent{\bf Abstract}:
we show that the logarithmic trace of Szeg\"o projectors introduced by K. Hirachi \cite{Hi04} for CR structures and extended in  \cite{lB05} to contact structures vanishes identically. }


%%%%%%%%%%%  Introduction - notations  %%%%%%%%%%%%%

\bigskip\bigskip In  \cite{Hi04} K. Hirachi showed that the logarithmic trace of the Szeg\"o projector is an invariant of the CR structure. In \cite{lB05} I showed that it is also defined for generalized Szeg\"o projectors associated to a contact structure  (definitions recalled below, sect.\ref{LBsz}), that it is a contact invariant, and that it vanishes if the base manifold is a 3-sphere, with arbitrary contact structure (not necessarily the canonical one). Here we show that it always vanishes. For this use the fact that this logarithmic trace is the residual trace of the identity (definitions recalled below, sect.\ref{LBtr}), and show that it always vanishes, because the Toeplitz algebra associated to a contact structure can be embedded in the Toeplitz algebra of a sphere, where the identity maps of all `good' Toeplitz modules have zero residual trace.

\section{Notations}
We first recall the notions that we will use. Most of the material below in \S1-5 is not new; we have just recalled briefly the definitions and useful properties, and send back to the literature for further details (cf. \cite{lH71,lH85,MS74,SKK73}).

If $X$ is a smooth manifold we denote $T^\bullet X\subset T^*X$ the set of non-zero covectors. A complex subspace $Z$ corresponds to an ideal $I_Z\subset C^\infty(T^\bullet X,\mathbb{C})$). $Z$ is conic (homogeneous) if it is generated by homogeneous functions. It is smooth if $I_Z$ is locally generated by $k=$ codim$\,Z$ functions with linearly independent derivatives. If $Z$ is smooth, it is involutive if $I_Z$  is stable by the Poisson bracket (in local coordinates $\{f,g\}=\sum \frac{\partial f}{\partial\xi_j}\frac{\partial g}{\partial x_j}-\frac{\partial f}{\partial x_j}\frac{\partial g}{\partial\xi_j}$); it is $\gg0$ if locally $I_Z$ has generators $u_i,v_j$ ($1\le i\le p,1\le j\le q,  p+q=$codim$\,Z$) such that the $v_j$ are real, $u_i$ complex, and the matrix $\frac1i(\{u_k,\bar{u}_l)$ is hermitian $\gg0$. The real part $Z_R$ is then a smooth real submanifold of $T^\bullet X$, whose ideal is generated by the $\rm{Re}\,u_i,\rm{Im}\,u_i,v_j$. If $\mathcal{I}$ is $\gg0$, it is exactly determined by its formal germ (Taylor expansion) along the set of real points of $Z_R$.

\medskip A Fourier integral operator (FIO) from $Y$ to $X$, is a linear operator from functions or distributions on $Y$ to same on $X$ defined as a locally finite sum of oscillating integrals 
$$
f\mapsto Ff(x)= \int e^{i\phi(x,y,\theta)} a(x,y,\theta) f(y) \ d\theta dy\;,
$$
where $\phi$ is a phase function (homogeneous w.r. to $\theta$), $a$ a symbol function. Here we will only consider regular symbols, i.e. which are asymptotic sums $a\sim\sum_{k\ge0} a_{m-k}$ where $a_{m-k}$ is homogeneous of degree $m-k$, $k$ an integer. 
$m$ could be any complex number. There is also a notion of vector FIO, acting on sections of vector bundles, which we will not use here.
 
The canonical relation of $F$ is the image of the critical set of $\phi\  (d_\theta\phi=0)$ by the differential map $(x,y,\theta)\mapsto (d_x\phi,-d_y\phi)$ - this is always assumed to be an immersion from the critical set onto a Lagrangian sub-manifold of $T^\bullet X\times T^\bullet Y^0$ (the sign $^0$ means that we have reversed the sign of the canonical symplectic form; likewise if $\mathcal{A}$ is a ring, $\mathcal{A^0}$ denotes the opposite ring). We will also use FIO with complex positive phase function: then the canonical relation is defined by its ideal (the set of complex functions $u(x,\xi,y,\eta$ which lie in the ideal generated by the coefficients of $d_\theta\phi,\xi-d_x\phi,\eta+d_y\phi$ which do not depend on $\theta$); it should not be confused with its set of real points.

Following H\"ormander \cite{lH71,lH85}, the degree of $F$ is defined as $\mbox{deg\,} F=\mbox{deg\,}(ad\theta)-\frac14
(n_x+n_y+2n_\theta)$, with $n_x,n_y,n_\theta$ the dimensions of the $x,y,\theta$-spaces, $\mbox{deg\,}(ad\theta)$ the degree of the differential form $ad\theta=\mbox{deg\,} a+\sum\nu_k$ ($\nu_k=\mbox{deg\,}\theta_k$, usually $\nu_k=1$ but could be any real number - not all $0$); this only depends on $F$ and not on its representation by oscillating integrals). In what follows we will always require that the degree be an integer (which implies $m\in \mathbb{Z}/4$).

%%%%%%%%%%% Adapted Fourier Integral Operators %%%%%%%%%%%%%

\section{Adapted Fourier Integral Operators}

The Toeplitz operators and Toeplitz algebras used here, associated to a CR or a contact structure, were introduced and studied in \cite{lB79.2,BG81}, using the analysis of the singularity of the Szeg\"o kernel (cf. \cite{lB76.5,mK77}), or in a weaker form, the ``Hermite calculus" of  \cite{lB74.4,vG74}. The terminology ``adapted" is taken from \cite{BG81}: lacking anything better I have kept it.

For $k=1,2$, let $\Sigma_k\subset T^\bullet X_k$ be smooth symplectic sub-cones of $T^\bullet X_k$ and $u:\Sigma_1\to\Sigma_2$ an isomorphism. 

\begin{definition} A Fourier integral operator $A$ is adapted to $u$ if its canonical relation $\mathcal{C}$  is complex $\gg0$, with real part the graph of $u$. It is elliptic if its principal symbol does not vanish (on $\mbox{Re\,}\mathcal{C}$).\end{definition}

As above, a conic complex Lagrangian sub-manifold $\Lambda$ of $T^\bullet X$ is $\gg0$ if its defining ideal $\mathcal{I}_\Lambda$ is locally generated by $n=\mbox{dim\,} X$ homogeneous functions $u_1,\dots,u_n (n=\mbox{dim\,} X)$ with independent derivatives, $u_1,\dots,u_k$ complex, $u_{k+1},\dots,u_n$ real for some $k$ ($1\le k\le n$), and the matrix $(\frac1i\{u_p, \bar{u}_q\})_{1\le p,q\le k}$ is hermitian $\gg0$; equivalently; the intersection $\mathcal{C}\cap\bar{\mathcal{C}}$ is clean, and on the tangent bundle the hermitian form $\frac1i\omega(U,\bar{V})$ is positive, with kernel the complexification of the tangent space of $\mbox{Re\,}\mathcal{C}$. 

Pseudo-differential operators are a special case of adapted FIO ($X_1=X_2= X, u=Id_{T^\bullet X}$);  so are Toeplitz operators on a contact manifold (see below).

 Adapted FIO always exist (cf \cite{BG81}), more precisely
\begin{proposition}\label{LBadapt} For any symplectic isomorphism  $u$ as above, there exists an elliptic FIO adapted to $u$.\end{proposition}
In fact if $\Lambda$ is a complex $\gg0$ Lagrangian sub-manifold of $T^\bullet X$, in particular if it is real, it can always be defined by a global phase function with positive imaginary part ($\rm{Im}\,\phi\succsim \mbox{dist\,}(.,\Lambda_R)^2)$ living on $T^\bullet X$: it is easy to see that such phase functions exist locally, and the positivity condition makes it possible to glue things together using a homogeneous partition of the unity. Once one has chosen a global phase function, it is obviously always possible to choose an elliptic symbol - of any prescribed degree (cf. also \cite{BG81})\footnote{the intrinsic differential-geometric description of the symbol is elaborate: it is a section of a line bundle whose definition incorporates half densities and  the Maslov index or an elaboration of this in the case of complex canonical relations. However on real manifolds this line bundle is always topologically trivial}. Note that elliptic only means that the top symbol is invertible on the real part of the canonical relation, not that the operator is invertible mod. smoothing operators (for this the canonical relation must be real:  $X_1,X_2$ have the same dimension, $\Sigma_k=T^\bullet X_k$ and $\mathcal{C}$ is the graph of an isomorphism).

%%%%%%%%%%%%%%%    model     %%%%%%%%%%%%%%%%%%%%
\section{Model Example}
\medskip Here is a generic example of adapted FIO: let $X_1,X_2,Z$ be three vector spaces 
$$\Sigma_k=T_{X_k}X_k\times T^\bullet Z\subset T^\bullet(X_k\times Z)\quad(k=1,2, \ T_{X_k}X_k \mbox{ the zero section}),$$
\begin{equation}
u\  \mbox{the identity map Id}_{T^\bullet Z}:\ \Sigma_1\to\Sigma_2 \;.
\end{equation}
If $\mathcal{C}\gg0$ is a complex canonical relation with  real part the graph of $u$, the complex formal germ along $\Sigma$ of the restriction to $\mathcal{C}$ of the projection $(x,\xi,z,\zeta,z',\zeta',y,\eta)\mapsto (x,z,\zeta',y)$ is an isomorphism (the dimensions are right, and it is an immersion: if $v=(0,\xi,0,\zeta,z',0,0,\eta)$ is a complex vector with zero projection, it is orthogonal to $\bar{v}$ (because these vectors form a real Lagrangian space), so if it is tangent to $\mathcal{C\gg0}$, it is tangent to the real part, i.e. the diagonal of $T^\bullet Z$, and this obviously implies $v=0$).

So we can choose the phase function as 
$$
\phi=<z-z',\zeta'> + iq(x,z,\zeta',y)
$$
where $q$ is smooth complex function of $x,z,\zeta',y$ alone, homogeneous of degree $1$ w.r. to $\zeta'$, vanishing of order $2$ for $x=y=0$, and $\rm{Re\,}q\succsim(x^2+y^2)\; |\zeta'|$ (it is easy to check that conversely any such phase function corresponds to a positive adapted canonical relation as above). The operator is
\begin{equation}
Ff(x,z)=\int e^{i<z-z',\zeta'>-q(x,z,\zeta',y)} a(x,z,\zeta',y) f(z',y) d\zeta'dz'dy \;,
\end{equation}
with $a$ a symbol as above.

Since any symplectic sub-cone of a cotangent manifold is always locally equivalent to $T^\bullet Z\subset T^\bullet(X\times Z)$, the model above is universal i.e. any adapted FIO is micro-locally equivalent to $F_1\circ A\circ F_2$ where $F_1,F_2$ are elliptic invertible FIO with real canonical relations, graphs of local symplectic isomorphisms, and $F$ is as the model above.

For adapted FIO the H\"ormander degree coincides with the degree in the scale of Sobolev spaces, i.e. if $F$ is of degree $s$ it is continuous $H^m(Y)\to H^{m-\rm{Re}\,s}(X)$; this is easily seen on the model example above ($F$ is $L^2$ continuous if its degree is $0$ i.e. $a$ is of degree $-\frac14 (n_x+n_y)$). (This not true for general FIOs - in fact for a FIO with a real canonical relation $\mathcal{C}$, this is only true if $\mathcal{C}$ is locally the graph of a symplectic isomorphism.)

\bigskip The following result also immediately follows from the positivity condition:
\begin{proposition}  Let $X_1,X_2,X_3$ be three manifolds, $\Sigma_k\subset T^\bullet X_k$ symplectic sub-cones, $u$ resp. $v$ a homogeneous symplectic isomorphism $X_1\to X_2$ resp. $X_2\to X_3$, $F,G$  FIO (with compact support) adapted to $u,v$. Then $G\circ F$ is adapted to $v\circ u$; its canonical relation is transversally defined and positive. It is elliptic if $F$ and $G$ are elliptic.\end{proposition}

This is mentioned in \cite{BG81}; the crux of the matter is that if $Q(y)$ is a quadratic form with $\gg0$ real part, the integral $\int e^{-|\xi| Q(y)}dy$ does not vanish: it is an elliptic symbol of degree $-\frac12 n_y$, equal to $\rm{disc}(\frac Q\pi)^{-\frac12} |\xi|^{-\frac12 n_y}$.

%%%%%%%%%%%%%%%%   Generalized Szeg\"o projectors   %%%%%%%%%
\section{Generalized Szeg\"o projectors\label{LBsz}}

These were called ``Toeplitz projectors" in \cite{lB02,lB05}. C. Epstein suggested the present name, which is better. References: {\cite{lB76.5,lB79.2,BG81,lB97.3}}.

\begin{definition} Let $X$ be a manifold, $\Sigma\subset T^\bullet X$ a symplectic sub-cone. A generalized Szeg\"o projector (associated to $\Sigma$) is an elliptic FIO $S$ adapted to Id$_\Sigma$ which is a projector ($S^2=S$)\end{definition}

(Note that ``elliptic'' (or ``of degree $0$'') is part of the definition; otherwise there exist many non-elliptic projectors, of degree $>0$ as FIOs). The case we are most interested in is the case where $\Sigma$ is the half line bundle corresponding to a contact structure on $X$ (i.e. the set of positive multiples of the contact form). But everything works as well in the slightly more general setting above.

We will not require here that $S$ be an orthogonal projector; this makes sense anyway only once one has chosen a smooth density to define $L^2$-norms.

If $S$ is a generalized Szeg\"o projector, its canonical relation $\mathcal{C}\subset T^\bullet X\times T^\bullet X$ is idempotent, positive, and can be described as follows: the first projection is a complex positive involutive manifold $Z_+$ with real part  $\Sigma$; the second projection is a complex negative manifold $Z_-$ with real part $\Sigma$ ($Z_-=\bar{Z_+}$ if $S$ is selfadjoint). The characteristic foliations define fibrations $Z_\pm\to \Sigma$ (the fibers are the characteristic leaves;  they have each only one real point so are ``contractible" (they vanish immediately in imaginary domain), and there is no topological problem for them to build a fibration). Finally we have $\mathcal{C}=Z_+\times_\Sigma Z_-$
 
 \medskip  Generalized Szeg\"o projectors always exist, so as orthogonal ones (cf. \cite{BG81,lB97.3}). 
As mentioned in \cite{lB05}, generalized Szeg\"o projectors mod. smoothing operators form a soft sheaf on $\Sigma$, i.e. any such projector defined near a closed conic subset of $T^\bullet X$ or $\Sigma$ is the restriction of a globally defined such projector.

%%%%%%%%%%%%%%%%   trace  %%%%%% %%%%%%%%%
\section{Residual trace and logarithmic trace \label{LBtr}}

The residual trace was introduced by M. Wodzicki \cite{mW87}. It was extended to Toeplitz operators and suitable Fourier integral operators by V. Guillemin \cite{vG93} (cf. also \cite{vG96,sZ99}). It is related to the first example of `exotic' trace given by J. Dixmier \cite{jD66}.

Let $\mathcal{C}$ be a canonical relation in $T^\bullet X\times T^\bullet X$. A family $A_s (s\in\mathbb{C})$ of FIOs of degree $s$ belonging to  $\mathcal{C}$ is holomorphic if $s\mapsto (\Delta)^{-\frac s2} A_s$ is a holomorphic map from $\mathcal{C}$ to FIO of fixed degree (in the obvious sense). If $A_s$ has compact support, the trace $\mbox{tr\,} A_S$ is then well defined and depends holomorphically on $s$ if $\mbox{Re\,} s$ is small enough ($A_s$ is then of trace class). Often, e.g. if the canonical relation is real analytic, this will extend as a meromorphic function of $s$ on the whole complex plane, but this is not very easy to use because the poles are hard to locate and usually not simple poles.

\begin{proposition}  If $\mathcal{C}$ is adapted to the identity Id$_\Sigma$, with $\Sigma\subset T^\bullet X$  a symplectic sub-cone, and $A_s$ a holomorphic family, as above, then $\mbox{tr\,} A_s$ has at most simple poles are the points $s=-n-\mbox{deg\,} A_0+k$,  $k\ge0$ an integer, $n=\frac12\mbox{dim\,} \Sigma$ (the degree $\mbox{deg\,} A_0$ is defined as above)\end{proposition}

Proof: this is obviously true is $A_s$ is of degree $-\infty$ (there is no pole at all). In general we can write $As$ as a sum of FIO with small micro-support  (mod smoothing operators), and a canonical transformation reduces us to the model case, where result is immediate.

\begin{definition} If $A$ is a FIO adapted to $\mathcal{C}$, the residual trace $\mbox{tr}_{res} A$ is the residue at $s=0$ of any holomorphic family $A_s$ as above, with $A=A_0$. \end{definition}

This does not depend on the choice of a family $A_s$: indeed if $A_0=0$, the family $A_s$ is divisible by $s$ i.e. $A_s=sB_s$ where $B_s$ is another holomorphic family, and since $\mbox{tr\,} B_s$ has only simple poles, $\mbox{tr\,} A_s$ has no pole at all at $s=0$.

\begin{proposition} The residual trace is a trace, i.e. if $A$ and $B$ are adapted Fourier integral operators, we have $\mbox{tr}_{res} AB = \mbox{tr}_{res} BA$. \end{proposition}

Indeed with the notations above $\mbox{tr\,}AB_s$ and $\mbox{tr\,}B_sA$ are well defined and equal for $\mbox{Re\,} s$ small, so their meromorphic extensions and poles coincide.

%%%%%   Logtrace  %%%%%%%%%
\bigskip\noindent{\bf Logarithmic trace (contact case)}

\medskip Let $\Sigma\subset T^\bullet X$ be a symplectic half-line bundle, defining a contact structure on $X$. A complex canonical relation $\mathcal{C}\gg0$ adapted to $\rm{Id}_\Sigma$ is always the conormal bundle of a complex hypersurface $Y$ of $X\times X$, with real part the diagonal (rather the positive half)\footnote{indeed a complex vertical vector $v$ is as before orthogonal to $\bar{v}$; if it is tangent to $\mathcal{C}$ at a point of $\mathcal{C}_R=\rm{diag}\,\Sigma$, it is tangent to the real part $\mathcal{C}_R=\rm{diag}\,\Sigma$ since $\mathcal{C}\gg0$, but this implies that $v$ is the radial vector, because the radial vector is the only vertical vector tangent to $\Sigma$. So the projection $\mathcal{C}\to X\times X$ is of maximal rank $2n-1$ and the image is a hypersurface.}, so if $A$ is a FIO adapted to $\rm{Id}_\Sigma$, its Schwartz kernel can be defined by a one dimensional Fourier integral:
$$
\tilde{A}(x,y)=\int_0^\infty e^{-T\phi(x,y)}a(x,y,T)\;dT\,,
$$
with $\phi=0$ an equation of the hypersurface $Y$, $\phi=0$ on the diagonal, $\rm{Re}\,\phi \ge cst\;\rm{dist}(.,\rm{diag})^2$, and $a$ is a symbol: $a\sim\sum_{k\le N}a_k(x,y)T^{k-1}$  ($N=\rm{deg}\,A$).

Its singularity has a typically holonomic form:
\begin{equation}
f(x,y)(\phi+0)^{-N}+g(x,y)\; \mbox{Log\,}(\frac1{\phi+0}) \;,
\end{equation}
with $f,g$ smooth functions on $X\times X$, and in particular $g(x,x)=a_0(x,x)$.

\begin{proposition} 
With notations as above, the residual trace of $A$ is  the trace of the logarithmic coefficient:
\begin{equation}
\mbox{tr}_{res} A = \int_X g(x,x)\; .
\end{equation}
\end{proposition}
An obvious holomorphic family extending $A$ (mod. a smoothing operator) is the family $A_s$ with Schwartz kernel
$$
\tilde{A_s}(x,y)=\int_1^\infty e^{-T\phi(x,y)}a(x,y,T)\,T^s dT \;.
$$
Since $a_s(x,x,T)\sim\sum_{k\le N} T^{s+k-1}a_k(x,x)$ and  $\phi(x,x)=0$, we get 
$$
\tilde{A}_s(x,x)\sim\sum\frac{a_k(x,x)}{s+k}\;,
$$
with an obvious notation: the meromorphic extension of the trace just has just simple poles at each integer $j\ge-N$, with residue $\int_X a_{-j}(x,x)$. In particular the residue for $s=0$ is the logarithmic trace. 

\medskip The residual trace is also equal to the logarithmic trace in the case of pseudo-differential operators, or in the model case. In general the residual trace is well defined, but  I do not know if the logarithmic coefficient can be reasonably defined e.g. if the projection $\Sigma\to X$ is not of constant rank. For the equality with the residual trace, and for theorem \ref{LBc=1} below, the sign is important: the logarithmic trace is the integral of the coefficient of $\rm{Log\,}\frac1\phi$, not the opposite.

%%%%%   Ttrace  %%%%%%%%%
%\bigskip\noinden
\section{\bf Trace on a Toeplitz algebra $\mathcal{A}$ and on $\mbox{End\,}_{\mathcal{A}}(M)$}
If $S$ is a generalized Szeg\"o projector associated to $\Sigma\subset T^\bullet X$. The corresponding Toeplitz operators are the Fourier integral operators of the form $T_P=SPS$, $P$ a pseudo-differential operator (equivalently, the set of Fourier integral operators $A$ with the same canonical relation, such that $A=SAS$). They form an algebra $\mathcal{A}$ on which the residual trace is a trace: $\mbox{tr}_{res} AB = \mbox{tr}_{res} BA$. Mod. smoothing operators, this can be localized, and the Toeplitz algebra $\mathcal{A}_\Sigma$ is this quotient; it is a sheaf on $\Sigma$ (or rather on the basis). 
\begin{proposition}\label{LBindep} 
The Toeplitz algebra, so as the residual trace of $S$ only depend on $\Sigma$ and not on the embedding $\Sigma_k\subset T^\bullet X$
\end{proposition}
Indeed if $\Sigma\to T^\bullet X'$ is another embedding, $S'$ a corresponding Szeg\"o projector,  it follows from prop. \ref{LBadapt} that there exist elliptic adapted FIO $F,F'$ from $X$ to $X'$ resp. $X'$ to $X$ such that $F=S'FS,F'=SF'S', FF'\sim S,F'F\sim S$ so $S,S'$ have the same residual trace, and $A\mapsto FAF'$ is an isomorphism of the two Toeplitz algebras. In the lemma wa could as well embed $\Sigma$ in another symplectic cone endowed with a Toeplitz structure.

\medskip The definition of the residual trace extends in an obvious way to $\mbox{End\,}_\mathcal{A}(M)$ when $M$ is a free $\mathcal{A}$-module ($\mbox{End\,}_{\mathcal{A}}(M)$ is isomorphic to a matrix algebra with coefficients in $\mathcal{A}^0$, where $\mbox{tr}_{res}(a_{ij}) = \sum \mbox{tr}_{res} a_{ii}$ is obviously a trace, independent of the choice of a basis of $M$). It extends further to the case where $M$ is locally free (a direct summand of a free module), and to the case where $M$ admits (locally or globally) a finite locally free resolution: if
$$
0\to L^N \xrightarrow{\text d} \dots L^1\xrightarrow{\text d} L^0 \to 0
$$
is such a resolution, i.e. a complex of are locally free $\mathcal{A}$-modules $L^j$, exact in degree $\neq 0$, with a given isomorphism $\epsilon : L^0/dL^1 \to M$:
Then $\mbox{End\,}_\mathcal{A}(M)$ is isomorphic to ${\rm Rhom}^0(L,L)$, i.e. any $a\in\mbox{End\,}_\mathcal{A}(M)$ extends as morphism 
$\tilde{a}$ of complexes of $L$ ($\tilde{a}=(a_j), a_j\in\mbox{End\,}_\mathcal{A}(L^j)$) (if $a$ has compact support, $\tilde{a}$ can be chosen with compact support). Any two such extensions $\tilde{a},\tilde{a}'$ differ by a super-commutator $[d,s]=ds+sd$ ($s$, of degree $1$, can be chosen with compact support if $\tilde{a}-\tilde{a}'$ has compact support). The super-trace 
$$
\mbox{supertr}_{res}(\tilde{a}) = \sum (-1)^j\mbox{tr}_{res}(a_j)
$$
is then well defined; it only depends on $a$, because the super-trace of a super-commutator $[s,d]$ vanishes: this defines the trace in $\mbox{End\,}(M)$. Below we will use ``good" modules, i.e. which have a global finite locally free resolution for which this is already seen on the principal symbols (i.e. $M$ and the $L_j$ are equipped with good filtrations for which $\mbox{gr\,} d$ is a locally free resolution of $\mbox{gr\,} M$ (in the analytic setting this always exists if $M$ is coherent and has a global good filtration, and the base manifold is projective).

%%%%%   Alt-trace  %%%%%%%%%
\bigskip\noindent{\bf Alternative description of the residual trace}

\medskip Let $S$ be a generalized Szeg\"o projector associated to a symplectic cone $\Sigma\subset T^\bullet X$. 
Then the left annihilator of $S$ is a $\gg0$  ideal $\mathcal{I}$ in the pseudo-differential algebra of $X$; its characteristic set is the first projection $\Sigma_+$ of the complex canonical relation of $S$; as mentioned earlier it is involutive $\gg0$ with real part $\Sigma$.
(there is a symmetric statement for the right annihilator).

\begin{proposition}\label{LBM} Let $M$ be the $\mathcal{E}_X$-module $M = \mathcal{E}/\mathcal{I}$. Then the Toeplitz algebra is canonically isomorphic to $\mbox{End\,}_\mathcal{E}(M)$. \end{proposition}

Proof: let  $e_M$ be the image of $1\in\mathcal{E}$ (it is a generator of $M$). It is elementary that $\mbox{End\,}_\mathcal{E}(M)$ is identified with opposite algebra $([\mathcal{E}:\mathcal{I}]/\mathcal{I})^0$ where $[\mathcal{E}:\mathcal{I}]$ denotes the set of $P\in\mathcal{E}$ such that $\mathcal{I} P\subset \mathcal{I}$ (to $P$ corresponds the endomorphism $a_P$ such that $a_P(e)=Pe$).
It is also immediate that the map $u$ which to $a_P$ assigns the Toeplitz operator $T_P=SPS=PS$ is an isomorphism (both algebras have a complete filtration by degrees, and the associated graded algebra in both cases is the algebra of symbols on $\Sigma$); clearly $u$ is a homomorphism of algebras, of degree $\le 0$, and $\mbox{gr\,} u = {\rm Id}$.
 
Now $M$ is certainly ``good'' in the sense above: it is locally defined by transverse equations and has, locally, a resolution whose symbol is a Koszul complex. So the residual trace is well defined on $\mbox{End\,}_\mathcal{E}(M)$. Since the trace on an algebra of Toeplitz type is unique up to a constant factor, there exists a constant $C$ such that
\begin{equation}
\mbox{tr}_{res} a_P = C \,\mbox{tr}_{res} T_P \;.
\end{equation}
Below we only  only need $C\neq0$; however with the conventions above we have:
\begin{theorem}\label{LBc=1}
The constant $C$ above is equal to one ($C=1$).
\end{theorem}
Proof: to the resolution of $M$ above corresponds a complex of pseudo-differential operators 
\begin{equation} \label{LBresol}
0 \to C^\infty(x)\xrightarrow{\text D} C^\infty(X,E_1)\to \dots \xrightarrow{\text D} C^\infty(X,E_N) \to 0 \;,
\end{equation}
exact in degree $>0$ and whose homology in degree $0$ is the range of $S$ (mod. smoothing operators), i.e. there exists a micro-local operator $E$ on $(E_k)$ such that $DE+ED \sim 1-S$ (cf. \cite{lB74.4}; $E$ is a pseudo-differential operator of type $\frac12$, not a ``classical" pseudo-differential operator, but it preserves micro-supports).

It is elementary that one can modify D, E, and if need be $S$, by smoothing operators so that 
\eqref{LBresol} is exact (on global sections) in degree $\neq0$, and $\ker D_0$ is the range of $S$. Then if $a_s=(a_k,s)$ is a holomorphic family of pseudo-differential homomorphisms of degree $s$, and $T_{a_s}$ is the Toeplitz operator $T_{a_s}=a_{0,s}|_{\ker S}$, we have $\mbox{tr\,} T_{a_s} = \sum (-1)^k\mbox{tr\,} a_{k,s}$ for $\mbox{Re\,} s\ll0$ hence also equality for the meromorphic extensions and residues.

\medskip Here is an alternative proof (slightly more in the spirit of the paper because it really uses operators  mod. $C^\infty$ rather than true operators). Notice first that theorem \ref{LBc=1} is (micro) local, and since locally all bundles are trivial, we can reason by induction on $\rm{codim}\,\Sigma$. Thus it is enough to check the formula for one example, where $\mbox{codim\,}\Sigma=2$. Note also that above we had embedded in the algebra of pseudo-differential operators on a manifold, but we could just as well embed in another Toeplitz algebra. 

We choose $\Sigma$ corresponding to the standard contact (CR) sphere $S_{2n-1}$ of $\mathbb{C}^n$, embedded as the diameter $z_1=0$ in the sphere of $\mathbb{C}^{n+1}$ 

The Toeplitz space $H_{n+1}$ is the space of holomorphic functions in the unit ball of $\mathbb{C}^{n+1}$ (more correctly: their restrictions to the sphere); we choose $H_n$ the subspace of functions independent of $z_1$. There is an obvious resolution:  $0\to H_{n+1}\to H_{n+1}\to 0$ ($H_n=\ker \partial_1$). We choose on $H_n$ the operator $a$, restriction of $\rho^-n$, with $\rho=\sum z_j\partial_{z_j}$ (this is the simplest operator with a nonzero residual - our convention is that $\rho^{\sigma}$ kills constant functions for all $\sigma$).

\begin{lemma} On the sphere $S_{2n-1}$ the residual trace of the Toeplitz operator $\rho^{-n}$ is \ 
$\mbox{tr}_{res} \rho^{-n}=\frac1{(n-1)!}$.
\end{lemma}
Proof: the standard Szeg\"o kernel is $S(z,w)= \frac1{\rm{vol\,}S_{2n-1}} (1-z.\bar{w})^{-n}$. Now we have the obvious identity
$$
\rho(\rho+1)\dots(\rho+n-1)\ \mbox{Log\,} \frac1{1-z.\bar{w}} = (n-1)! ((1-z\bar{w})^{-n}-1) \;,
$$
so that the leading coefficient of the logarithmic part of $\rho^{-n}$ is $\frac 1{\rm{vol\,}S_{2n-1}(n-1)!}$, whose integral over the sphere is  $\frac1{(n-1)!}$.

\medskip On the sphere $S_{2n+1}$ we have $\partial_1\rho = (\rho+1) \partial_1$ so we choose for $(\tilde{a})$ the pair $(a_0=\rho^{-n},a_1=(\rho+1)^{-n}$). Since terms of degree $<-n-1$ do not contribute on $S_{2n+1}$, the super-trace is
$$
\mbox{supertr}_{res}(\tilde{a}) = \mbox{tr}_{res} (\rho^{-n}-(\rho+1)^{-n})=\mbox{tr}_{res}(n \rho^{-n-1}) = \frac n{n!}=\mbox{tr}_{res} a \;.
$$

%%%%%%%%%%%  Embedding  %%%%%%%%%%%%%
\section{Embedding}

If $\Sigma_1,\Sigma_2$ are two symplectic cones, with contact basis $X_1,X_2$, symplectic embeddings $\Sigma_1\to\Sigma_2$ exactly correspond to contact embeddings $X_1\to X_2$, i.e. an embedding $u:X_1\to X_2$ such that the inverse image $u^*(\lambda_2)$ is a positive multiple of the contact form $\lambda_1$ (the corresponding symplectic map take the section $u^*\lambda_2$ of $T^\bullet X_1$ to the section $\lambda_2$ of $T^\bullet X_2$. With this in mind we have

\begin{lemma}\label{LBembed} If $X$ is a compact oriented contact manifold, it can be embedded in the standard contact sphere.\end{lemma}

Proof: the standard contact $(2N-1)$-sphere of radius $R$ has coordinates $x_j,y_j\ (1\le j\le N$, $\sum x_j^2+y_j^2=R^2$) and contact form $\lambda=\sum x_jdy_j-y_jdx_j$ (or a positive multiple of this).
If $X$ is a compact contact manifold, its contact form can always be written $2\sum_2^mx_j dy_j$ for some suitable choice of smooth functions $x_j,y_j$, or just as well $\sum_1^m x_jdy_j-y_jdx_j$, setting for instance $x_1=1,y_1=\sum_2^m x_jy_j$ ($m$ may be larger than the dimension). Adding suitably many other pairs $(x_j,y_j)$ with $y_j=0,x_N=(R^2-\sum x_j^2+y_j^2)^{\frac12}$, for $R$ large enough, we get an embedding in a contact sphere of radius $R$.

\begin{theorem} For any generalized Szeg\"o projector $\Sigma$ associated to a symplectic cone with compact basis, we have $\mbox{tr}_{res} S=0$. In particular if $\Sigma$ corresponds to a contact structure on a compact manifold, the logarithmic trace vanishes.
\end{theorem}

By lemma \ref{LBembed} we can suppose that $\Sigma$ is embedded in the symplectic cone of a standard odd contact sphere. Let $\mathcal{B}$ be the canonical Toeplitz algebra on the sphere: then by prop. \ref{LBM}, the Toeplitz algebra $\mathcal{A}$ of $S$ is isomorphic to $\mbox{End}_\mathcal{B}(M)$ where $M$ is a suitable good $\mathcal{B}$ module. Now on the sphere any good locally free $\mathcal{B}$-module is free (any complex vector bundle on an odd sphere is trivial), and the Szeg\"o projector has no logarithmic term, so  and $\mbox{tr}_{res} 1_M=0$, for any free hence also for any good $\mathcal{B}$-module $M$.

\bigskip This result is rather negative since it means that the logarithmic trace cannot define new invariants distinguishing CR or contact manifolds. Note however that that it is not completely trivial: it holds for the Toeplitz algebras associated to a CR or contact structure, as constructed in \cite{BG81}, but a contact manifold carries many other star algebras which are locally isomorphic to the Toeplitz algebra (I showed in \cite{lB02} how Fedosov's classification of star products \cite{bF94} can be adapted to classify these algebras). Any such algebra $\mathcal{A}$ carries a canonical trace, because the residual trace is invariant by all isomorphisms, so that local traces glue together. If the contact basis is compact, the trace $\mbox{tr}_{res} 1_\mathcal{A}$ is well defined, but there are easy examples showing that it is not always zero.

%%%%%%%%%%%  bibli  %%%%%%%%%%%%%

\bigskip{\sc Universit\'e Pierre et Marie Curie - Paris 6, Analyse Alg\'ebrique, Institut de Math. de Jussieu, Case 82 - 4, place Jussieu, 75252 Paris Cedex 05, France}

{\it e-mail:  boutet@math.Jussieu.fr}

\end{document}